\documentclass[12pt]{article}

\normalbaselines
\catcode `\@=11
\@addtoreset{equation}{section}

\def\section{\@startsection {section}{1}{\z@}{-3.5ex plus -1ex minus
     -.2ex}{2.3ex plus .2ex}{\normalsize\bf}}
\def\subsection{\@startsection{subsection}{2}{\z@}{-3.25ex plus -1ex minus
 -.2ex}{1.5ex plus .2ex}{\normalsize\bf}}
\def\thebibliography#1{\section*{References\markboth
  {REFERENCES}{REFERENCES}}\list
  {[\arabic{enumi}]}{\settowidth\labelwidth{[#1]}\leftmargin\labelwidth
  \advance\leftmargin\labelsep
  \usecounter{enumi}}
  \def\newblock{\hskip .11em plus .33em minus -.07em}
  \sloppy
  \sfcode`\.=1000\relax}
 
\catcode `\@=12
\font\ab=msbm10 scaled\magstep1% for C, R, Z
\font\aba=msbm7 % for C, Z, R in small format
\newcommand{\gz}{\mbox{\ab Z}}%Z double
\newcommand{\gf}{\mbox{\ab F}}%F double
\newcommand{\gl}{\mbox{\ab L}}%L double
\newcommand{\gk}{\mbox{\ab K}}%K double
\newcommand{\gcp}{\mbox{\ab CP}}%CP double
\newcommand{\gpl}{\mbox{\ab PL}}%PL double
\newcommand{\gpk}{\mbox{\ab PK}}%PL double
\newcommand{\gc}{\mbox{\ab C}}%C double
\newcommand{\gcm}{\mbox{\aba C}}% C double for exponent
\newcommand{\Ka}{K\"ahler}
\newcommand{\Gras}{\mbox{$G_n({\gc}^{m+n})$}}% Grassmannian
\newcommand{\men}{\mbox{$\widetilde {\bf M}$}}% m tilde
\newtheorem{com}{Comment}

\newtheorem{pr}{Proposition}
\newtheorem{topic}{Topic}

\begin{document}
\vspace*{1.5cm}
\noindent
\begin{center}
{\bf COHERENT STATES, TRANSITION AMPLITUDES AND EMBEDDINGS}
\vspace{1.0cm}\\
\end{center}
\noindent

\vspace{1.3cm}

\noindent
%\hspace*{1in}
\begin{center}
\begin{minipage}{14cm}
S. Berceanu \\
Institute of Physics and Nuclear Engineering,\\
Department of Theoretical Physics,\\ P. O. Box MG-6, Bucharest-Magurele,
Romania;\\
E-mail: Berceanu@Theor1.IFA.Ro;  Berceanu@Roifa.IFA.Ro\\
\makebox[3mm]{ }\\
\end{minipage}
\end{center}

\vspace*{0.5cm}

\begin{abstract}
\noindent
%Type here the abstract of your contribution, if you want to have one.
The transition  amplitudes between  coherent states on the coherent
state manifold \men~ are expressed in terms of the embedding of \men~
into a projective Hilbert space \gpl .  Consequences for the dimension of
\gpl~ and a simple geometric interpretation of Calabi's diastasis follows.

\end{abstract}

\vspace{1cm}

\setcounter{section}{-1}
\section{\hspace{-4mm}.\hspace{2mm}{\bf Introduction and preliminaries}}

The coherent states \cite{klauder}
 offer a powerful framework  to formulate  a link between
classical and quantum mechanics \cite{sbaa,sbl}. Simultaneously,
 the coherent state approach
furnishes an appealing recipe \cite{berezin}
 for the geometric quantization \cite{kost}.
 However,
a "physical" motivation of the group theoretic generalisation 
  of Heisenberg-Weyl's group to arbitrary
Lie groups due to Perelomov \cite{per} is still missing. On the
other side, a simple geometric description of the coherent states, and firstly
of transition amplitudes and transition probabilities \cite{wigner,bargmann},
 is well-suited.
 The proposed way
to attain this goal is the embedding of the coherent state manifold into an
 adequate projective Hilbert space. The importance  of this embedding was
already emphasised \cite{raw,od}.  A very simple answer to these
questions is obtained
formulating the problem in the language of
 complex geometry \cite{chern}, fibre bundles \cite{hus}
 and algebraic geometry \cite{hirz}.
  
In this talk we shall be concerned with the following topics:
1) the geometric meaning of the transition amplitudes;
2) angles, distances and coherent states;
 3)  the geometric meaning of Calabi's diastasis;
4) Kodaira embedding and coherent states.
Elsewhere \cite{sbpol} we have been concerned with the questions:
 5) the relationship 
between geodesics and coherent states; 6) a geometric
 characterisation of the polar divisor, i.e. the set of coherent vectors
 orthogonal to a fixed vector.
Putting the answers to all these questions together, we get a better
 understanding \cite{sb} of the coherent
states. A full illustration of the problems 1)-6) in the case  of the complex
Grassmann manifold  \Gras~ is given in Ref. \cite{viitor}.

%\section{\hspace{-4mm}\hspace{2mm}{\bf Preliminaries}}
%\vspace*{1cm}
\subsection{\hspace{-5mm}.\hspace{2mm} The coherent states}

Let  $\pi$ be an  unitary irreducible representation,  $G$ a Lie group and
 \gk~  a Hilbert space. Let the orbit  $\widetilde {\bf M} =
{\widetilde
\pi (G)}\vert \widetilde {\psi_0>}$, where $\vert  \psi_0> \in \gk$ and
$\xi :\gk\rightarrow\gpk$ is the projection
$ \xi \vert  \psi>=\vert \widetilde\psi>$. Then there is
 the diffeomorphism $\widetilde {\bf M} \approx G/K $, where $K$ is the
 stationary group of  $\vert \widetilde\psi_0>$. If $\iota :
\widetilde {\bf M} \hookrightarrow \gpl$ is an embedding,
 then $\widetilde {\bf M}$ is called {\it coherent state manifold}.
 If $\vert  \psi_0> \equiv \vert j>$ is an extreme weight
vector, then for  compact connected simply connected Lie groups
$G$,
 $\widetilde {\bf M}$  is a K\"ahler manifold and the celebrated
 Borel-Weil-Bott \cite{bwb} theorem furnishes both the representation
 $\pi_j$ and the representation space $\gl=\gk^*_j$, where
$X^*$ denotes the dual of the vector space $X$.
 If a local section $\sigma:\widetilde {\bf
M}\rightarrow{\cal S}(\gk)$ in the unit sphere in \gk~ is constructed,
 then the    holomorphic line
bundle ${\bf M}' =\sigma(\widetilde {\bf M})$ is  associated by a
 holomorphic character  $\chi$ of the
parabolic subgroup $ P $  of the complexification  $G^{\gcm}$ of $G$.

 The   coherent vectors \cite{per}, which belong to the {\it coherent
 vector manifold} {\bf M} \cite{sbcag}, are introduced as
\begin{equation}
\vert
Z,j>=\exp\sum_{{\varphi}\in\Delta^+_n}(Z_{\varphi}F^+_{\varphi})\vert
j> ,~{\vert
\underline{Z}>}=< Z \vert Z>^{-1/2}\vert Z>\in {\bf M} ,\label{zzz}
\end{equation}
 where $\Delta^+_n$ are the positive non-compact roots,
  $ Z\equiv (Z_ \varphi )
\in \gc^n$ are the local coordinates in  neighbourhood
 of $|j>$ corresponding to $Z=0$
  and  $n$ is the dimension
of the manifold \men . We remember that
  $F^+_{\varphi}\vert j>\neq  0,~F^-_{\varphi}\vert j> = 0,~  \varphi
\in\Delta^+_n. $

  Below $<\underline{Z}' \vert\underline{Z}>$ denotes the   hermitian
 scalar
product of holomorphic sections in the line bundle ${\bf M}$ in different
points of the manifold \men .

\subsection{\hspace{-5mm}.\hspace{2mm} The manifold}

The first study of compact homogeneous complex manifolds was done by H. C.
Wang \cite{wang}, who completely classified  the $C-$spaces, i.e.
the simply connected
compact homogeneous manifolds. If $G$ is a connected semisimple Lie group, then
a K\"ahlerian $C-$space is
necessarily of the form $G/C(T)$, where $T$ is a toral subgroup of $G$
and $C(T)$ is the centralizer of $T$ in $G$. Then {\it every compact
 homogeneous K\"ahler
manifold is a K\"ahlerian direct product of a K\"ahlerian $C$-space and a flat
complex torus} (cf. Matsushima's theorem, see e.g. Note
24 pp. 373-375 in Ref. \cite{kn}).

 The following theorem summarises some properties of flag manifolds with
 significance for the present paper \cite{wolf}.

 {\it Let $X_c=G^{\gcm} /P$ be a complex manifold, where
$G^{\gcm}$ is a complex semisimple Lie group and $P$
 is a parabolic subgroup. The following conditions are equivalent:

a)  $X_c=G^{\gcm}/P$ is compact;

b) $X_c$  is a complex connected \Ka~ manifold;

c) $X_c$ is a projective variety;

d) $X_c$ is a closed $G^{\gcm}$ orbit in a projective representation;

e) $X_c$ is a Hodge manifold and all homogeneous Hodge manifolds
are of this type.}

We remember that the manifold \men~ is called a {\it Hodge} manifold
 (\Ka~ manifold of restricted type) if the \Ka~ two-form $\omega$ is integral, 
i.e. $\omega\in H^2(\men ,\gz)$.\newpage

\subsection{\hspace{-5mm}.\hspace{2mm} The embedding} 
\label{embedding}

A holomorphic line bundle ${\bf M}'$ on a compact complex manifold  \men~
is said {\it very ample} \cite{ss} if: the set of divisors is without
 base points, i.e. there exists
a finite set of global sections $s_1,\ldots ,s_N\in \Gamma (\men, \bf{M'})$
such that for each $m\in \men$ at least one $s_j(m)$ is not zero, and
 the holomorphic map $\iota_{{\bf M}'}:\men\hookrightarrow\gcp^{N-1}$ 
given by
\begin{equation}
  \iota_{{\bf M}'}=[s_1(m),\ldots ,s_N(m)]\label{scufund}
\end{equation}
is a holomorphic embedding. So, 
 $\iota_{{\bf M}'}:\men \hookrightarrow\gcp^{N-1}$ is an embedding
 if \cite{gh} :

\AA$_1)$ the set of divisors is without base points;

\AA$_2)$ the differential of $\iota$ is nowhere degenerate;

\AA$_3)$ $\iota$ is one-one, i.e. for any $m, m' \in \men$
there exists $s\in H^0(\men,\cal{O}({\bf M'}))$ such that $s(m)=0$
and $s(m')\not= 0$, where
 $\cal{O}$ denotes the sheaf of holomorphic sections.

The line bundle ${\bf M}'$ is said to be {\it ample} if there exists a positive
integer $r_0$ such that ${\bf M}'^r$ is very ample for all $r\geq r_0$. Note
that if ${\bf M'}$ is an ample line bundle on \men , then \men ~must be
 projective-algebraic by Chow's theorem, hence \men~ is \Ka .

The holomorphic line bundle ${\bf M'}$ is said to be {\it positive} if
on ${\bf M}'$ can be given a hermitian metric
 $ds^2\in C^{\infty}(\men,{\bf M'}^*\times \overline{\bf M}'^*)$ such
that $\sqrt{-1}\Theta$~ is positive, where $\Theta~ $is the curvature form of
 the hermitian connection. If in local coordinates the two-form
 $\omega\in \wedge ^{1,1}$ is $\omega=
 \sqrt{-1}\sum g_{ik}dz_i\wedge d\bar{z}_k$, then
  $\omega$  is positive if the matrix $[g_{ik}]$ is positive definite.

The concepts of ampleness and positivity for line bundles coincide. The
following theorem \cite{ss} summarises the properties of ample line bundles
 that are needed in this paper.

 {\it Let ${\bf M}'$ be a holomorphic line bundle  on a compact
complex manifold \men . The following conditions are equivalent:

a) ${\bf M}'$ is positive;

b) for all coherent analytic
 sheaves ${\cal S}$ on \men~ there exists a
 positive integer $m_0(\cal{S})$ such that
 $H^i(\men,{\cal S}\otimes {\bf M}'^m)=0$ for $i>0,~m\geq m_0(\cal{S})$
 (the vanishing theorem of Kodaira);

c) there exists a positive integer $m_0$ such that for all $m\geq m_0$,
there is an embedding  $\iota_{\bf M}:\men \hookrightarrow\gcp^{N-1}$ for
some $N\geq n$ such that ${\bf M}={\bf M'}^m$ is projectively induced, i.e.
${\bf M}=\iota ^*[1]$;

d) \men ~ is a Hodge manifold  (the embedding theorem of Kodaira);

e) in particular, if \men~ is also homogeneous, then  \men~ is a flag 
 manifold.}

In the condition of  case {\it e)}\/, i.e. when \men~ is a
 homogeneous K\"ahler manifold,   the exact description of the embedding
$\iota_{\bf M}:\men \hookrightarrow\gcp^{N-1}$ is
furnished by the  Borel-Weil-Bott theorem \cite{bwb}.
 The dimension of the
 representation is given by the Riemann-Roch-Hirzebruch theorem. The same
 result can be obtained using the coherent states, as will be  seen later in
 Proposition \ref{bigtm}.  Here $[1]$ denotes the hyperplane bundle.

Now we discuss the construction of the embedding $\iota :\men\hookrightarrow
\gpl$ for noncompact manifolds. Then the projective Hilbert space is infinite
dimensional \cite{kobi}.

Let  \gf~ be the Hilbert space of square integrable holomorphic $n-$forms
on \men . Then $\gl =\gf ^*$.  Let $z=(z_1,\ldots ,z_n)$ be a local
 coordinate system. Let $\iota '$ be the mapping which sends $z$ into an
 element $\iota '(z)$ of \gl~  defined by the paring $<\iota' (z),f>=
f^*(z)$, where $f(z_1,\ldots ,z_n)=f^*dz_1\wedge \cdots \wedge dz_n\wedge
d\overline{z}_1\wedge \cdots \wedge d\overline{z}_n$. Then $\iota '(z)\not=
0$ if a condition analogous to condition \AA$_1$) in the noncompact case
 is satisfied. Then $\iota
=\xi \circ \iota '$ is independent of local coordinates and is continuous and
complex analytic. 

If $K$ is the kernel $2n-$form on $\men \times \overline{\men}$ then the \Ka~
metric of Kobayashi \cite{kobi} is $ds^2=\sum\partial^2\log K^*/
\partial z_i\partial \overline{z}_jdz_id\overline{z}_j$, where
 $K(z,\overline{z})=K^*(z,\overline{z})dz_1\wedge \cdots\wedge dz_n\wedge
d\overline{z}_1\wedge \cdots\wedge d\overline{z}_n$. 

The analogous of conditions \AA$_1$)-\AA$_3$) used by Kobayashi in the
 noncompact  case are:

A$_1$) for any $z\in \men$, there exists a square integrable $n-$form
$f$ such that $f(z)\not= 0$;

A$_2$) for every holomorphic vector $Z$ at $z$ there exists a square integrable
$n-$form $f$ such that $f(z)=0$ and $Z(f^*)\not= 0$;

A$_3$) if $z$ and $z'$ are two distinct points of \men, then there is a
 $n-$form $f$ such that $f(z)=0$ and $f(z')\not= 0$.

Kobayashi has shown that condition A$_1$) implies $ds^2=\iota ^*(ds^2_{FS})$,
while A$_2$) and A$_3$) imply that $\iota$ is also an embedding.

Rawnsley \cite{raw} has globalized the definition of coherent states
 including also the non-homogeneous  \Ka~ manifolds. He has shown  that
 $\omega_{\men}-\iota ^*\omega_{FS}=
{\displaystyle\frac{1}{2\pi  i}}\bar{\partial}\partial\eta$, where
$ \omega_{FS}$  is the fundamental two-form on the complex projective
space. So, if
$\eta$ is harmonic, then $\iota$ is K\"ahlerian and an immersion. 
 For regular hermitian line bundle, in particular
for homogeneous K\"ahler manifolds and homogeneous quantization,
 $\eta $ is constant  and $\omega$ is
the pull-back of $\omega_{FS}$. For the
complex torus $T=\gc^n/\Gamma$, $T$ is Hodge if and only if the Riemann
conditions are satisfied \cite{chern}. The projectively induced line
bundles correspond to $\iota$ an embedding.

\section{\hspace{-4mm}.\hspace{2mm}  The geometric meaning of the transition 
amplitude}

\begin{topic}: find a geometric meaning of the transition probability on
coherent state manifold.\end{topic}
 \begin{pr}\label{tramp}
 Let $\vert
\underline{Z}>$ as in (\ref{zzz}), where Z parametrizes the coherent state
 manifold
in the ${\cal V}_0 \subset \widetilde {\bf M} $ and let us suppose that 
the coherent state manifold admits  the embedding
 $\iota:\widetilde {\bf M}\hookrightarrow \gpl$. Then the angle
\begin{equation} 
\theta \equiv \arccos \vert <\underline{Z'} \vert\underline{Z}>\vert
,\label{un1}
\end{equation}
is equal to the Cayley distance on the geodesic joining $\iota (Z'),\iota (Z),$
 where $Z', Z \in {\cal V}_0$,
\begin{equation}
\theta = d_c(\iota (Z'),\iota (Z)).\label{unul1}
\end{equation}
 
More generally, it is true the following relation (Cauchy formula)
 
$$<\underline{Z}' \vert\underline{Z}>=
\frac{(\iota (Z'),\iota (Z))}{\|\iota (Z')\|\|\iota (Z)\|}.$$
\end{pr} 
{\it Proof:} We discuss here the case of compact manifolds.
 The embedding (\ref{scufund})    is realised in the case of the
coherent state manifold \men~ by the formula
\begin{equation}
\label{scut}\iota (Z)=[|\underline{Z}>]~.
\end{equation}

Because  the manifold \men~ admits a embedding into the projective Hilbert
 space \gpl , the
line bundle ${\bf M}'$ is a positive one.
 The theorem from Section
\ref{embedding} is applied. It follows that there is a power $m_0$
of the positive line bundle ${\bf M}'$  such that the coherent vector manifold
verifies the relation ${\bf M}={\bf M}'^{m_0}$.  The holomorphic line bundle
 {\bf M} of coherent vectors is the
 pull-back $\iota^*$ of
the hyperplane bundle $[1]$ of \gpl, the dual bundle of the tautological
line bundle of $\gpk^*$, i.e. ${\bf M}=\iota^*[1]$. The analytic line 
bundle  ${\bf M}$ is {\it projectively induced} (see p. 139 in Ref.
\cite{hirz}).
 
 In the Proposition \ref{tramp}, (.,.) is the scalar product in \gk. 
If $\xi :\gk\backslash  \{0\}
\rightarrow
 \gpk,~\xi : \omega \rightarrow [\omega ]$, then the elliptic hermitian
 Cayley \cite{cay} distance is
\begin{equation}\label{cd} 
 d_c([\omega'],[\omega ])=\arccos{ \vert (\omega ',\omega )\vert \over  \| 
\omega '\| \| \omega \| }.
\end{equation}
 
 The noncompact case is treated similarly.

 For completeness, we remember here the notion of tautological 
line bundle \cite{hirz}.
  $[1]=[1_n]$ is the $\gc^{\star}$-bundle defined by
the cocycle $\{g_{ij}\}=\{z_jz_i^{-1}\}$, where 
$[z_0,\ldots ,z_n]$ are the homogeneous coordinates for the complex projective 
space $\gcp^n$. $\gcp^{n+1}\setminus \{0\}$ is a principal bundle
  with structure
group $\gc^{\star}$ which is associated to the $U(1)$-bundle $[1_n]^{-1}=
[-1_n]$.
The principal bundle $U(n+1)/U(n)$ over the Grassmann manifold
$G_1(\gc^{n+1})=\gcp^{n}=SU(n+1)/S(U(n)\times U(1))$ is associated to the
tautological (universal)   bundle over $\gcp^n$.

\begin{com} {\bf (The distances in Quantum Mechanics: variations on a theme
 by Cayley)}
~ \end{com} The  Cayley distance (\ref{cd}) has been used independently in
 Quantum
 Mechanics by many authors \cite{wick,fivel,aa}. The Cayley distance
(\ref{cd}) is useful in the
geodesic approach. The elliptic hermitian distance $d_c$ of two points
given by eq. (\ref{cd})
is one half the arc of the great circle connecting the corresponding points on
the Riemann sphere \cite{coo}. Some authors \cite{study} prefer instead of
 eq. (\ref{cd}) the definition
\begin{equation}\label{str} 
 d'_c([\omega'],[\omega ])=2\arccos{ \vert (\omega ',\omega )\vert \over  \|
\omega '\| \| \omega \| }.
\end{equation}

 The (Bargmann \cite{bargmann}) distance $d_b$, used by
 Prevost and Vall\' ee \cite{prv}
in the context of  coherent states,
 
$$d^2_b([\omega'],[\omega ])=2(1-\cos\,d_c([\omega'],[\omega ])),$$
is equivalent with $d_c :{2\sqrt {2}}/ \pi d_c \leq d_b\leq d_c.$

Defining the inner product $(\alpha\beta)$ of two rays as the absolute
value of the scalar product $<\alpha|\beta >$, a ``distance''
 $\rho_{\alpha\beta}$  between two rays is introduced by formula (2) at page
232 in Ref.  \cite{wick}:
\begin{equation}\label{wi}
\cos\, (\frac{1}{2}\rho_{\alpha\beta})=(\alpha\beta)=|<\alpha|\beta >|,~
0\leq\rho_{\alpha\beta}\leq\pi .
\end{equation}

The connection between geodesics in the space of rays and probability 
 transition  is commented in  \S \ ``Some remarks on ray space'' of Ref.
 \cite{wick}.  One shows that if $\alpha$ and $\beta$ are
not orthogonal ($\rho_{\alpha\beta}<\pi$) a condition for $\gamma$, stronger 
than linear dependence, is that $\gamma$ should lie on the geodesic arc
 connecting $\alpha$ to $\beta$ and in  this case
\begin{equation}
\rho_{\alpha\beta}=\rho_{\alpha\gamma}+\rho_{\gamma\beta}~.\label{sum}
\end{equation}

Formula (\ref{wi}) is identical with eq. (6) in Ref.  \cite{aa}   :
\begin{equation}
|<\psi|\phi >|^2=\cos ^2(\frac{1}{2}\theta )~,\label{ananahar}
\end{equation}
where $|\phi>,~|\psi >$ are points in $\gc^{N+1}$ and $\theta$
is the distance joining $\xi (|\psi >)$ and $\xi (|\phi >)$.

In fact, {\it formula (\ref{ananahar}) of  Anandan and Aharonov \cite{aa} and
respectively formula (\ref{wi})  of Wick \cite{wick} were knew from the last
 century} (see Ref.  \cite{cay} pp. 584, 590).
So, {\it formula (\ref{ananahar}) is nothing else than the definition
 (\ref{str}) of the distance on the projective space.}

\section{\hspace{-4mm}.\hspace{2mm}  Angles, distances and coherent states}
 
\begin{topic}: find those manifolds $\widetilde {\bf M} $ for which the angle
given by eq. (\ref{un1}) 
 is a distance on $\widetilde {\bf M} $.\end{topic}
 
\begin{pr} Let $\widetilde {\bf M} $ be a coherent state manifold 
parametrized as in
(\ref{zzz}). Then the angle given by eq. (\ref{un1}) 
is a distance on $\widetilde {\bf M} $ iff $\widetilde {\bf M} $ is a symmetric
space of rank 1.
 \end{pr}
{\it Proof:} The problem is reduced to that of two-point homogeneous
spaces, which are known \cite{wolfi}.
 
\begin{com}Generally, the distance $\delta$ on a manifold is
 greater than the angle $\theta$ defined by eq. (\ref{un1}),
 $\delta \geq \theta $, but 
infinitesimally, $d\delta=d\theta .$ 
\end{com}

\section{\hspace{-4mm}.\hspace{2mm} The geometric meaning of Calabi's
 diastasis}

\begin{topic}: find a geometric meaning of Calabi's diastasis \cite{cal},
 used by
Cahen, Gutt,
Rawnsley \cite{cgr} in the context of coherent states,
 $D(Z',Z)=-2\log\, \vert <\underline{Z}' \vert\underline{Z}>\vert .$
 \end{topic}
\begin{pr}
 The diastasis distance $D(Z',Z)$ between $Z', Z \in {\cal V}
_0
 \subset \widetilde {\bf M}$ is related to the geodesic distance
 $\theta = d_c(\iota (Z'),\iota (Z))$,
 where $\iota :\widetilde {\bf M} \hookrightarrow \gpl,$ by
 
\begin{equation}
D(Z',Z)= -2\log\, \cos\,\theta .
\end{equation}
 
 If $\widetilde {\bf M}_n$ is noncompact and $\iota ':\widetilde {\bf M}_n 
 \hookrightarrow \gcp^{N-1,1}=
  SU(N,1)/S(U(N)\times U(1))$ ~
${\mbox (\iota :\widetilde {\bf M}_n \hookrightarrow \gpl )}$
and $\delta _n (\theta
  _n)$ is the length of the geodesic joining $\iota '(Z'),\iota '(Z)~$
 (resp. $\iota (Z'),\iota (Z))$,
then $$\cos\,\theta _n=
  (\cosh\,\delta _n)^{-1}=e^{-D/2}.$$
\end{pr} 
  {\it Proof:} The proposition is a direct consequence of
 Proposition \ref{tramp}.

\begin{com} The remark \cite{sbpol} {\bf the polar divisor = cut locus}
for manifolds \men~ of symmetric type gives a geometric description of the
 domain of definition of  Calabi's diastasis.
\end{com}

\section{\hspace{-4mm}.\hspace{2mm}  Kodaira embedding and coherent states}
 
\begin{topic}: characterise the relationship of the smallest number $N$
in the Kodaira embedding $\iota :\widetilde {\bf M} \hookrightarrow
 \gcp^{N-1}$ and
 the compact complex manifold $\widetilde {\bf M}$.
 
\end{topic}
 
\begin{pr}\label{bigtm} For coherent state manifolds
 $\widetilde {\bf M} \approx  G/K $ which have a flag
 manifold structure, the following 7 numbers are equal:
 
 1) the maximal number of orthogonal coherent vectors on $\widetilde {\bf M}$;
 
 2) the number of holomorphic global sections in the holomorphic line bundle
{\bf M} with base $\widetilde {\bf M}$;
 
 3) the dimension of the fundamental  representation in the Borel-Weil-Bott
 theorem;

 4) the minimal $N$ appearing in the Kodaira embedding theorem,
 $\iota :\widetilde {\bf M} \hookrightarrow  \gcp^{N-1}$;
 
 5) the number of critical points of the energy function $f_H$ attached to a 
 Hamiltonian $H$
 linear in the generators of the Cartan algebra of G, with unequal
 coefficients;
 
 6) the Euler-Poincar\'e characteristic of $\widetilde {\bf M} \approx G/K,~ 
\chi (\widetilde {\bf M})=
 [W_G]/[W_H]$, where $[W_G]= card\,W_G$, and $ W_G$ is the Weyl group of G;
 
 7) the number of Borel-Morse cells which appear in the CW-complex 
 decomposition of $\widetilde {\bf M}$.
\end{pr}
 
 {\it Proof:} Use theorems 1, 2 in Ref.  \cite{sbcag}
 where it is proved that $f_H$
 is a perfect Morse function, the Cauchy formula and the Borel-Weil-Bott
theorem \cite{bwb}. Remark that $\chi (G/K)>0$ iff
 ${\rm Rank}\,G={\rm Rank}\,K$, cf.
 to a classical result of Hopf and Samelson \cite{hs}.
 
 \begin{com}
 The Weil prequantization condition is the condition to have a 
 Kodaira embedding,
 i.e. the algebraic manifold to be Hodge.
\end{com}

 The author expresses his thanks to the Organising Committee
 to invite him to the  XIV Workshop on Geometrical
 Methods in Physics at  Bia\l owie\.za.  Discussions during the Workshop,
 especially with
Professors  M. Cahen,  J. Klauder, W. Lisiecki, M. Schlichenmaier and T.
Wurzbacher  are acknowledged.
The constant interest of Professor L. Boutet de Monvel is kindly acknowledged.

\end{document}